\documentclass[a4paper]{article}
\usepackage{amsmath,amsthm,amsfonts}

\newtheorem{thm}{Theorem}[section]
\newtheorem{lemma}[thm]{Lemma}
\newtheorem{prop}[thm]{Proposition}

\newcommand{\prlabel}[1]{\label{#1}}%\fbox{#1}}
\newcommand{\prbibitem}[2]{\bibitem[#1]{#2}}%\fbox{#2}}

\begin{document}

\title{On the asymptotic geometry of area--preserving maps}

\author{Leonid Polterovich\\School of Mathematical Sciences\\Tel Aviv
  University\\69978 Tel Aviv\\Israel \and Karl Friedrich Siburg\\Fakult\"at
  f\"ur Mathematik\\Ruhr--Universit\"at Bochum\\44780 Bochum\\Germany}

\date{\textsc{Preliminary Version}}

\maketitle

\section{Statement and discussion of the results}

Let $M$ be an open connected oriented 2--manifold endowed with an area form
$\omega$. We assume that the total area of $M$ with respect to $\omega$ is
infinite, i.e. $\int_M \omega= \infty$. Consider the group Ham$_c(M,\omega)$
of Hamiltonian diffeomorphisms of $M$ consisting of all time--1--maps of
time--periodic compactly supported Hamiltonians $H: \mathbb{S}^1\times M\to
\mathbb{R}$.  We write $\varphi_H^t$ for the Hamiltonian flow generated by
$H$.

We are interested in the asymptotic behaviour of one--parameter subgroups of
Ham$_c(M,\omega)$ with respect to Hofer's metric $d$ where \footnote{ We refer
  the reader to \cite{book}, \cite{intro} and \cite{con} for an introduction
  to Hofer's geometry.}
\[ d(\text{id},\varphi)= \inf\{ \int_0^1 
\max H_t-\min H_t\,dt\mid H \text{ has compact support and } \varphi_H^1=
\varphi\}. \]
Let $\mathcal{A}$ denote the Lie algebra of Ham$_c(M,\omega)$; it consists of
all compactly supported time--independent Hamiltonians on $M$.  Given $H\in
\mathcal{A}$, we are interested in the growth of the function $r_H:
[0,\infty)\to [0,\infty)$ defined by
\[ r_H(t)= d(\text{id},\varphi_H^t) \;. \]
By the triangle inequality for Hofer's metric we know that $r_H$ is
subadditive, i.e.
\[ r_H(t+s)\leq r_H(t)+ r_H(s) \;. \]
Therefore the limit
\[ \mu(H)= \lim_{t\to\infty} \frac{r_H(t)}{t} \]
is well defined. This quantity---which is called the asymptotic
non--minimality of the subgroup generated by $H$---was introduced in
\cite{geod}; see also \cite{con} for further discussion.

In particular, we see that $r_H$ grows at most linearly in $t$. On the other
hand, if $(M,\omega)$ is the standard Euclidean plane then a theorem by
Sikorav (\cite{sikorav,book}; see also Proposition \ref{jcstrick} below)
states that $r_H$ is bounded by a constant, and this constant depends only on
the diameter of the support supp($H$) of $H$.

In the present note, we show that for open surfaces of infinite area the
function $r_H$ is either bounded or behaves asymptotically linear. In order to
formulate our main result we need the following notion. Recall that a subset
$Z \subset M$ is called contractible in $M$ if the inclusion $Z
\hookrightarrow M$ is homotopic to the constant map which sends the whole $Z$
to a point in $M$.

\begin{thm}[Dichotomy Theorem]
  \prlabel{dicho} Let $(M,\omega)$ be an open surface of infinite area and $H:
  M\to \mathbb{R}$ a compactly supported Hamiltonian. Then the following
  dichotomy holds:
\begin{itemize}
\item If $\{H\neq 0\}$ is contractible in $M$ then the function $r_H$ is
  bounded; in particular, $\mu(H)=0$.
\item If $\{H\neq 0\}$ is not contractible in $M$ then the function $r_H$
  grows asymptotically linear, i.e., $\mu(H)>0$.
\end{itemize}
\end{thm}

In fact, it is even possible to calculate the precise value of $\mu(H)$ as the
difference of two distinguished critical values of $H$. In particular, one
gets examples of one--parameter subgroups whose asymptotic non--minimality
lies strictly between 0 and $\max H - \min H$ and can be calculated precisely.
As far as we know, this is the first series of examples of this type.

Assume that $M \neq \mathbb{R}^2$. Then $\pi_1(M)$ is nontrivial. Let
$\mathcal{L}$ be the set of all embedded non--contractible circles in $M$.
Then we define
\begin{align*}
  c_+(H)&= \sup_{L\in\mathcal{L}} \min_{x\in L} H(x)\\
  c_-(H)&= \inf_{L\in\mathcal{L}} \max_{x\in L} H(x)
\end{align*}
Since $H$ is compactly supported, it follows that $c_+(H) \geq 0$ and $c_-(H)
\leq 0$ (see Proposition \ref{sig}). Moreover, one can easily check (see
Proposition \ref{crit}) that $c_+(H)$ and $c_-(H)$ are critical values of $H$.

\begin{thm}
  \prlabel{value} The following equality holds:
\[ \mu(H)= c_+(H)- c_-(H) \;. \]
\end{thm}

We give a short outline of the proof of the two theorems. If $\{H\neq 0\}$ is
contractible then a version of the abovementioned theorem by Sikorav shows
that $r_H$ is bounded by a constant depending only on the ``size'' of
supp($H$). Thus we get the first statement of Theorem~\ref{dicho} (see Section
3 below). An elementary argument (see Proposition \ref{cont}) shows that
$c_-(H) = c_+(H) = 0$ if and only if $\{H \neq 0\}$ is contractible in $M$.
Thus, the second statement of Theorem~\ref{dicho} follows from
Theorem~\ref{value}. Theorem~\ref{value} consists of two parts. The inequality
$\mu(H)\geq c_+(H)- c_-(H)$ follows from a Lagrangian intersection result as
in \cite{imrn} (see Section 4). The proof of the reversed inequality uses a
trick, namely a decomposition of $\varphi_H^t$ into two commuting flows:
\[ \varphi_H^t= \Phi^t\circ \Psi^t= \Psi^t\circ \Phi^t \]
where $\Phi^t$ has contractible support. Hence, Sikorav's theorem yields that
$\Phi^t$ can be neglected in the calculation of $\mu(H)$ (see Section 3). Let
us emphasize that in order to apply Sikorav's argument we need that $M$ has
infinite area.

Intuitively, the two distinguished critical values $c_{\pm}(H)$ correspond to
the first homotopically nontrivial separatrices of $H$. More precisely,
$c_+(H)$ is the infimum of energy values $E$ such that the superlevel set
$H^{-1}([E,\infty))$ is contractible, and $c_-(H)$ is the supremum of $E$ with
contractible sublevel set $H^{-1}((-\infty,E])$. Consequently, the asymptotic
geometric behaviour of $\varphi_H^t$ depends only on the topology of the level
sets of $H$.

When $M$ is the cylinder, the lower bound on the asymptotic non--minimality
$\mu(H)$ in terms of the energy levels which carry non--contractible circles
was known \cite[9.B]{con}. Theorem~\ref{value} above shows that this bound is
sharp!

Concerning the function $r_H(t)$, we obtain the following picture. As long as
there are no non--constant periodic solutions we have
\[ r_H(t)= (\max H- \min H)\, t \;, \]
see \cite[II, Cor.\ 1.10]{geometry}. For large $t$, however,
\[ r_H(t)\sim (c_+(H)- c_-(H))\, t \;. \]
Consequently, if $\max H > c_+(H)$ or $c_-(H) > \min H$ there must be a
``phase transition'' in the behaviour of $r_H(t)$ from small $t$ to large $t$.

Let us conclude with a couple of open problems. First of all, it is not clear
if the Dichotomy Theorem holds true for surfaces of finite area, or, even more
ambituous, for higher--dimensional symplectic manifolds. We also do not know
how to deal with cyclic subgroups of Ham$_c(M,\omega)$, consisting of
time--1--maps of time--dependent Hamiltonians periodic in time. The reason is
the lack of an integral of motion which is essential for our arguments.
Finally, it would be interesting to have a dynamical interpretation of the
abovementioned change in the behaviour of $r_H(t)$.

\textit{Acknowledgement}: This work was done while the second author was
visiting Tel Aviv University in December 1998. He thanks the Minerva Center
for Geometry for the financial support.

\section{Some properties of $c_{\pm}(H)$}

In this section we sum up some useful elementary properties of $c_{\pm}(H)$
and deduce the second statement of Theorem~\ref{dicho} from
Theorem~\ref{value}. We start with some auxiliar facts and notions from
topology of open surfaces.

\begin{prop}
  \prlabel{aux1} An open subset $Z \subset M$ is contractible in $M$ if and
  only if every embedded circle which lies in $Z$ is contractible in $M$.
\end{prop}

\begin{proof}
  Assume without loss of generality that $Z$ is connected.  Since an open
  surface is a $K(\pi,1)$--space, the inclusion $Z \hookrightarrow M$ is
  homotopic to a point if and only if the homomorphism $\pi_1(Z) \to \pi_1(M)$
  is trivial \cite[8.1.11]{sp}. But $Z$ is also an open surface. Representing
  $Z$ as the union of an increasing chain of compact surfaces with boundary we
  see that there exists a countable system of embedded circles which generates
  $\pi_1(Z)$. Thus we get the proposition.
\end{proof}

Every contractible embedded circle $\gamma$ on $M$ bounds a unique closed disc
which we denote $D(\gamma)$. For a subset $X \subset M$ set \footnote{We write
  $\text{cl}(Z)$ and $\text{int}(Z)$ for the closure and the interior of a
  subset $Z$, respectively.}
\[ \text{hull}(X) = \text{cl}(\cup_{\gamma} D(\gamma)) \]
where $\gamma$ runs over all contractible embedded circles which are contained
in $X$.

\begin{prop}
  \prlabel{aux2} Let $N \subset M$ be a 2--dimensional compact submanifold
  with boundary. Assume that $\partial N$ is contractible in $M$. Then
  hull($N$) is the union of a finite number of pairwise disjoint closed discs
  whose boundaries are connected components of $\partial N$. In particular,
  $N$ is contractible in $M$.
\end{prop}

\begin{proof}
  It suffices to show that there exist pairwise disjoint closed embedded discs
  $D_1,\dots, D_k \subset M$ such that $\partial D_j$ is a component of
  $\partial N$ and $ N \subset \cup_j D_j$. Since $N$ is compact it has only
  finitely many connected components which we denote by $N_i$; let
  $\gamma_{ij}$ denote the boundary components of $N_i$. Now fix some $i$.
  Since $M$ is open and $\partial N$ contractible, at least one of the discs
  $D(\gamma_{ij})$ intersects the interior of $N_i$. Denote this disc by
  $D_i$.
  
  We claim that $D_i$ contains $N_i$. Indeed, pick any point $x \in
  \text{int}(D_i \cap N_i)$, and assume on the contrary that there exists a
  point $y \in \text{int}(N_i)\setminus D_i$. Since $N_i$ is connected, there
  exists a path $\alpha$ in the interior of $N_i$ which joins $x$ and $y$.
  But, since $x$ lies inside $D_i$ and $y$ outside $D_i$, $\alpha$ must
  intersect $\partial D_i \subset \partial N$ which is impossible. This
  contradiction proves the claim.
  
  Notice that for $i \neq j$ either $D_i$ and $D_j$ are disjoint, or one
  contains the other. So choose from the set $\{D_1,\ldots, D_p\}$ those discs
  which are maximal with respect to inclusion. This family of discs clearly
  satisfies all the requirements above.
\end{proof}

As an immediate corollary of the proposition let us mention that the hull of a
compact subset is compact. Indeed, each compact subset is contained in some
compact submanifold with boundary.

Let us return now to the quantities $c_{\pm}$. As before we assume that $M
\neq \mathbb{R}^2$. The next proposition is quite standard in the calculus of
variations and could be formulated, for instance, in the setting of the
Minimax Lemma in \cite[Sect.\ 3.2]{book}. For the convenience of the reader we
give a slightly more direct proof here.

\begin{prop}
  \prlabel{crit} $c_{\pm}(H)$ are critical values of $H$.
\end{prop}

\begin{proof}
  Assume on the contrary that $c_+(H)$ is a regular value of $H$. Then there
  exists a segment $[E_1,E_2]$ which consists of regular values of $H$ and
  such that $ E_1 < c_+(H) < E_2$. By definition of $c_+(H)$, the set $\{H
  \geq E_1\}$ contains a non--contractible circle. Since the gradient flow of
  $H$ takes $\{H \geq E_1\}$ into $\{H \geq E_2\}$ we conclude that $\{H \geq
  E_2 \}$ contains a non--contractible circle, too. Hence $c_+(H) \geq E_2$,
  in contradiction to the choice of $E_2$. The proof for $c_-(H)$ is
  analogous.
\end{proof}

\begin{prop}
  \prlabel{sig} $c_+(H) \geq 0$ and $c_-(H) \leq 0 $.
\end{prop}

\begin{proof}
  It suffices to show that $M\setminus {\rm supp}(H)$ contains a
  non--contractible curve from $\mathcal{L}$. Write $M= \cup_i N_i$, where
  $N_i$ are compact surfaces with boundary such that supp$(H)\subset
  \text{int}(N_1)$ and $N_i \subset \text{int}(N_{i+1})$ for all $i\geq 1$.
  
  If some boundary component of some $N_i$ is non--contractible in $M$ we are
  done. Assume therefore that all of them are contractible. Then Proposition
  \ref{aux2} implies that all $N_i$ are contractible in $M$. We conclude that
  $\pi_1(M) = 0$, i.e.\ $M = \mathbb{R}^2$, in contradiction to our standing
  assumption.
\end{proof}

\begin{prop}
  \prlabel{cont} $c_-(H) = c_+(H) = 0$ if and only if $\{H \neq 0\}$ is
  contractible in M.
\end{prop}

\begin{proof}
  We are going to apply Proposition \ref{aux1} with $Z= Z_H = \{H \neq 0\}$.
  
  If $Z_H$ is not contractible in $M$ then it contains a curve from
  $\mathcal{L}$. This curve lies either in $\{H >0\}$, which implies $c_+(H)
  >0$, or in $\{H < 0\}$, in which case $c_-(H) < 0$.
  
  Suppose now that $Z_H$ is contractible in $M$. Then $Z_H$ cannot contain a
  curve from $\mathcal{L}$. This means that every curve from $\mathcal{L}$
  intersects the set $\{H = 0\}$, so $c_+(H) \leq 0$ and $c_-(H) \geq 0$. But,
  as we have seen in Proposition \ref{sig}, $c_+(H)$ is non--negative and
  $c_-(H)$ is non--positive. Therefore $c_-(H) = c_+(H) = 0$.
\end{proof}

As a consequence we see that the second statement of Theorem~\ref{dicho}
follows from Theorem~\ref{value}.

\section{Decomposing Hamiltonian flows}

In this section we prove the following result.

\begin{thm}
  \prlabel{leq} The following inequality holds:
\[ \mu(H)\leq c_+(H)- c_-(H) \;. \]
Moreover, if $c_-(H) = c_+(H) = 0$ then $r_H$ is bounded.
\end{thm}

Together with Propositions \ref{sig} and \ref{cont} this implies the first
statement of Theorem~\ref{dicho}.

An essential ingredient of our approach is the following version of Sikorav's
theorem \cite{sikorav}; see also \cite[Sect.\ 5.6]{book}.

\begin{prop}
  \prlabel{jcstrick} Let $X \subset M $ be a finite union of pairwise disjoint
  closed discs, and let $F \in \mathcal{A}$ be a Hamiltonian function on $M$
  whose support is contained in the interior of $X$. Then
\[ d(\text{id}, \varphi_F^t) \leq 16\: \text{area}(X) \]
for every $t$.
\end{prop}

\begin{proof}
  When $M = (\mathbb{R}^2, dp \wedge dq)$ this is proved in
  \cite{sikorav,book}.  The case of a general open surface of infinite area
  can be reduced to this one as follows. Assume without loss of generality
  that $X\subset (M,\omega)$ consists of just one disc of area $A$. Let $D
  \subset (\mathbb{R}^2, dp\wedge dq)$ be the closed standard disc of area
  $A$. Since $M$ has infinite area, it is an easy consequence of the
  Dacorogna--Moser theorem (\cite{daco}, see also \cite[Sect.~1.6]{book}) that
  there exists a symplectic embedding
\[ i: (\mathbb{R}^2, dp \wedge dq) \hookrightarrow (M,\omega) \]
such that $i(D) = X$. Clearly, $i$ induces the natural homomorphism
\[ i_* : \text{Ham}_c(\mathbb{R}^2, dp \wedge dq) \to \text{Ham}_c(M, \omega)
\;. \]
It is important to notice that $i_*$ does not increase the corresponding Hofer
distances. Our flow $\varphi_F^t$ lies in the image of $i_*$, i.e.,
$\varphi_F^t = i_* (f_t)$ where $f_t$ is a one--parameter subgroup of
Ham$_c(\mathbb{R}^2, dp \wedge dq)$ whose Hamiltonian is supported in
int($D$). Thus, the desired inequality follows from Sikorav's original theorem
since $d(\text{id},\varphi_F^t) \leq d(\text{id},f_t) \leq 16\, A$.
\end{proof}

\begin{proof}[Proof of Theorem \ref{leq}]
  Let us decompose the flow $\varphi_H^t$ into two commuting flows as follows.
  Fix any $\epsilon >0$, and choose a smooth function $\rho:\mathbb{R}\to
  \mathbb{R}$ satisfying the following properties:
\begin{enumerate}
\item $\rho(s)=s$ if $c_-(H)-\epsilon \leq s \leq c_+(H)+\epsilon$
\item $\rho(s)=c_+(H)+2\epsilon$ if $s\geq c_+(H)+3\epsilon$
\item $\rho(s)=c_-(H)-2\epsilon$ if $s\leq c_-(H)-3\epsilon$
\item $0<\rho'(s)<1$ if $c_-(H)-3\epsilon < s < c_-(H)-\epsilon$ or
  $c_+(H)+\epsilon < s < c_+(H)+3\epsilon$
\end{enumerate}
Define the new Hamiltonians $K= \rho\circ H$ and $H_0= H-K$, and denote their
flows by $\Psi^t$ and $\Phi^t$, respectively. Then
\begin{equation}
\prlabel{commute}
\varphi_H^t= \Phi^t\circ \Psi^t= \Psi^t\circ \Phi^t \;.
\end{equation}
Observe that supp($H_0$) is contained in the set
\[ Z(\epsilon) = H^{-1} ((-\infty, c_-(H)-\epsilon] \cup 
[c_+(H) + \epsilon, \infty)) \;. \]

Pick any regular value $\kappa \in (0,\epsilon)$ of $H$. Then $Z(\kappa)$ is a
compact 2--dimensional submanifold with \textit{contractible} boundary. Denote
by $X$ the hull of $Z(\kappa)$. Proposition \ref{aux2} implies that $X$ is a
finite union of pairwise disjoint closed discs. Moreover, $Z(\kappa)$ is a
subset of supp($H$), so $X$ is contained in the hull of supp($H$).  Recall
that this hull is compact. Combining this with Proposition \ref{jcstrick}
above, we conclude that there is a constant $C>0$, depending only on supp($H$)
but not on $\epsilon$, such that
\[ d(\text{id},\Phi^t)\leq C \]
for every $t$.

On the other hand, $\Psi^t$ is generated by $K$ with $\max K- \min K= c_+(H)-
c_-(H)+ 4\epsilon$, hence
\[ d(\text{id},\Psi^t)\leq t\,(c_+(H)- c_-(H)+ 4\epsilon) \;. \]

Now, the relation \eqref{commute} implies that
\[ d(\text{id},\varphi_H^t)
\leq d(\text{id},\Phi^t)+
d(\text{id},\Psi^t) \;. \]
Therefore
\[ d(\text{id},\varphi_H^t)\leq C+ t\,(c_+(H)- c_-(H) + 4\epsilon) \]
for every $\epsilon > 0$, and the inequality in Theorem \ref{leq} follows.
Moreover, if $c_-(H) = c_+(H) = 0$ we get that $r_H(t)\leq C$ for all $t\geq
0$, and Theorem \ref{leq} is proven.
\end{proof}

\section{A lower bound on $\mu$ via Lagrangian intersections}

Recall from the introduction that $\mu(H)= \lim_{t\to\infty}
d(\text{id},\varphi_H^t)/t$. In the present section we prove the following
result.

\begin{thm}
  \prlabel{geq} We have the following inequality:
\[ \mu(H)\geq c_+(H)- c_-(H) \;. \]
\end{thm}

Together with Theorem \ref{leq} this completes the proof of
Theorem~\ref{value}, and thus that of Theorem~\ref{dicho}.

We will make use of the Lagrangian suspension construction and Lagrangian
intersection theory, similar to what is done in \cite{imrn}. It is convenient
to split Hofer's original definition for $d$ into two parts separating the
maximum and minimum. Define
\begin{align*}
  d_+(\text{id},\varphi)&= \inf_F\{\int_0^1 \max_x F_t\,dt\mid \varphi_F^1=
  \varphi\}\\
  d_-(\text{id},\varphi)&= \inf_F\{\int_0^1 -\min_x F_t\,dt\mid \varphi_F^1=
  \varphi\}
\end{align*}
where $F: \mathbb{S}^1\times M\to \mathbb{R}$ runs over all compactly
supported Hamiltonians generating $\varphi$. Then
\[ d(\text{id},\varphi)\geq d_+(\text{id},\varphi)+ d_-(\text{id},\varphi)
\;. \]
Recall that we consider $\varphi= \varphi_H^1$ where $H$ is autonomous.

\begin{lemma}
\prlabel{norm}
\begin{align*}
  d_+(\text{id},\varphi)&= \inf_F\{ \max_{t,x} F\mid \varphi_F^1= \varphi\}\\
  &= \inf_G\{ \max_{t,x} (H-G)\mid \varphi_G^1= \text{id}\}\\
  d_-(\text{id},\varphi)&= \inf_F\{ -\min_{t,x} F\mid \varphi_F^1= \varphi\}\\
  &= \inf_G\{ -\min_{t,x} (H-G)\mid \varphi_G^1= \text{id}\}
\end{align*}
\end{lemma}

\begin{proof}
  The first equalities are proved in \cite[\S 7]{imrn}, and the second ones in
  \cite[Lemma 3.A]{imrn}.
\end{proof}

\begin{lemma}
  \prlabel{loop} Suppose that $G: \mathbb{S}^1 \times M \to \mathbb{R}$ is a
  compactly supported Hamiltonian which generates the identity: $\varphi_G^1
  =$ id. Let $L \subset M$ be an embedded non--contractible circle. Then there
  exist $x_0 \in L$ and $t_0 \in \mathbb{S}^1$ such that $G(t_0,x_0) = 0$.
\end{lemma}

\begin{proof}[Proof of Theorem~\ref{geq}]
Since $\varphi_H^t = \varphi_{tH}^1$, it suffices to show that
$$d(\text{id},\varphi_H^1) \geq c_+(H) - c_-(H).$$
Fix an arbitrary $\epsilon
> 0$. Choose $L$ to be a non--contractible circle on $M$ such that $H|_L \geq
c_+(H)-\epsilon$. Lemmata \ref{norm} and \ref{loop} imply that
$d_+(\text{id},\varphi_H^1) \geq c_+(H) - \epsilon$.  Analogously,
$d_-(\text{id},\varphi_H^1) \geq -c_-(H) - \epsilon$. Thus
$d(\text{id},\varphi_H^1) \geq c_+(H) -c_-(H) - 2\epsilon$, for every
$\epsilon > 0$. Thus we get the desired inequality.
\end{proof}

\begin{proof}[Proof of Lemma~\ref{loop}]
  The proof goes along the lines of \cite{imrn}, and we only give a sketch
  here. The argument is devided into three steps.
  
  1) Choose a compact connected submanifold with boundary $N \subset M$ whose
  interior contains both $L$ and $\cup_t \text{supp}(G(.,t))$.  Let us perform
  the following surgery on $(M,\omega)$. We remove the complement to $N$ and
  attach to each boundary component of $N$ a cylindrical end of infinite area.
  Note that the loop of Hamiltonian diffeomorphisms $\varphi_G^t$ extends to
  this new surface. Therefore we can assume from the very beginning that
  $(M,\omega)$ has a finite number of ends and each end has infinite area.
  Such a surface $(M,\omega)$ is geometrically bounded (or tame) in the sense
  of Gromov's theory of pseudo--holomorphic curves (see \cite{hc}). This will
  enable us to apply Floer theory in Step 3 below.
  
  2) We claim that
\[ \pi_1({\rm Ham}_c (M,\omega)) = 0 \;. \]
This fact is well known to experts, however, as far as we know, no reference
is available.  Here is a sketch of the argument. Denote by Diff$_{c,0}(M)$,
respectively Symp$_{c,0}(M,\omega)$, the identity component of the group of
compactly supported diffeomorphisms, respectively symplectomorphisms, of $M$.
Consider the sequence
\[ \pi_1(\text{Ham}_c(M,\omega)) \to \pi_1(\text{Diff}_{c,0}(M)) \to
\pi_1(\text{Symp}_{c,0}(M,\omega)) \;. \]

The left arrow is a monomorphism (see \cite[Cor.\ 10.18(iii)]{intro} adjusted
to the non-compact case along the lines mentioned in the book). The right
arrow is an isomorphism; this follows from Moser's deformation argument with
parameters (cf.\ \cite[Sect.\ 3.2]{intro}). But it is shown in \cite{contr}
that $\pi_1(\text{Diff}_{c,0}(M)) = 0$. This completes the proof sketch of the
claim.

3) Finally, recall the so--called Lagrangian suspension construction for
Lagrangian submanifolds $L$ in a symplectic manifold $(M,\omega)$. Given $G:
\mathbb{S}^1\times M\to \mathbb{R}$ with $\varphi_G^1= \text{id}$, we consider
the embedding
\begin{align*}
  L\times \mathbb{S}^1&\to M\times T^*\mathbb{S}^1\\
  (x,t)&\mapsto (\varphi_G^t(x),t,-G(t,\varphi_G^t(x)))
\end{align*}
If we equip $M\times T^*\mathbb{S}^1$ with the split symplectic form
$\omega\oplus dr\wedge dt$ then the above map is a Lagrangian embedding. In
our case $L$ is a circle, and the image of the embedding is a Lagrangian torus
which we denote by $\mathcal{T}(G)$.

In view of Step 2 we know that the loop $\varphi_G^t, 0\leq t\leq 1$, is
homotopic to the constant loop at the identity. Hence the Lagrangian torus
$\mathcal{T}(G)$ is exact Lagrangian isotopic to $\mathcal{T}(0)=
\{(x,t,0)\mid x\in L, t\in \mathbb{S}^1 \}$. Moreover, since $L$ is
non--contractible, $\pi_2(M \times T^*S^1, \mathcal{T}(0)) = 0$. Then Floer
theory \cite{floer} guarantees the existence of an intersection point in
$\mathcal{T}(G)\cap \mathcal{T}(0)$, i.e., there are $x_0\in L$ and $t_0\in
\mathbb{S}^1$ such that
\[ G(t_0,x_0) = 0 \;. \]
This completes the proof of the lemma and finishes the proof of Theorem
\ref{geq}.
\end{proof}

\end{document}